\documentclass[12pt]{article}
\usepackage{amsmath}
\begin{document}
\title{Theorems on residues obtained by the division of powers\footnote{Delivered to the Berlin Academy
on February 13, 1755. Originally published as
\emph{Theoremata circa residua ex divisione potestatum relicta}, Novi Commentarii academiae scientiarum
Petropolitanae \textbf{7} (1761),
49-82, and reprinted in \emph{Leonhard Euler, Opera Omnia}, Series 1:
Opera mathematica, Volume 2, Birkh\"auser, 1992.
A copy of the original text is available
electronically at the Euler Archive, at www.eulerarchive.org. This paper
is E262 in the Enestr\"om index.}}
\author{Leonhard Euler\footnote{Date of translation: August 18, 2006.
Translated from the Latin by Jordan
Bell, 4th year undergraduate in Honours Mathematics, School of Mathematics
and Statistics, Carleton University, Ottawa, Ontario, Canada. Email:
jbell3@connect.carleton.ca.
This translation was written
during an NSERC USRA.}}
\date{}

\maketitle

\begin{center}
{\Large Theorem 1.}
\end{center}

1. If $p$ is a prime number and $a$ is prime to $p$, then
no term of the
geometric
progression $1,a,a^2,a^3,a^4,a^5,a^6$, etc. is divisible by the number $p$.

\begin{center}
{\Large Demonstration.}
\end{center}

This follows from Book VII, Prop. 26 of Euclid, where it is demonstrated
that if two numbers $a$ and $b$ are prime to $p$ then too that their product
$ab$ is prime to $p$. Thus with $a$ prime to $p$, by putting $b=a$, the
square $a^2$ will be prime to $p$; and then in turn by putting $b=a^2$;
likewise $b=a^3$, etc. Therefore no power of $a$ will be divisible by the prime
number $p$.

\begin{center}
{\Large Corollary 1.}
\end{center}

2. Therefore if each of the terms of the geometric progression
\[
1; a^2; a^3; a^4; a^5; a^6; a^7; a^8; \textrm{etc.}
\]
were divided by a prime number $p$, division never happens without
a residue, but rather a residue arises from each term.

\begin{center}
{\Large Scholion.}
\end{center}

3. I have resolved to carefully study 
the residues which emerge from the division of all
the terms of the given geometric progression by the prime number $p$.
First indeed, it is apparent from the nature of division that each of
these residues will be less than the number $p$; also, no residue
will be $=0$, because no term is divisible by $p$. For if
residues are produced which are greater than the number $p$, the way
in which they may be reduced to less than it is clear from arithmetic. 
Thus the residue of $p+r$ equals the residue $r$, and in general
the residue $np+r$ returns the residue $r$; and if $r$ is greater
than $p$, this residue is reprised as $r-p$, or $r-2p$, or $r-3p$, etc., until a
number less than $p$ is reached. Therefore all the residues 
$r \pm np$ reflect the same residue as $r$. In particular, all 
the residues may be referred to as 
positive numbers less than $p$.
In fact however, on many occasions it is convenient for negative residues
to be considered: as if $r$ were a residue remaining from division 
by a certain number $p$ such that $r<p$, the residue will be $r-p$, namely
a negative number; thus the residue of the positive
$r$ equals the residue of the negative $r-p$. 
In this way the residues can be exhibited such that none
exceed half of $p$, and on the other hand if a residue $r$ is 
larger than $\frac{1}{2}p$, its place is taken by the residue of the negative
$r-p$, since this will be less than half of $p$.

\begin{center}
{\Large Corollary 2.}
\end{center}

4. Because all the residues are integral numbers which are less than $p$, it
follows
that no more than $p-1$ different residues can arise. Since the terms geometric
series $1;a;a^2;a^3;a^4;a^5$; etc. consist of infinitely many different numbers,
it is necessary that multiple terms exhibit the same residues.

\begin{center}
{\Large Corollary 3.}
\end{center}

5. If $a^\mu$ and $a^\nu$ are two terms which produce the same residue $r$,
so that $a^\mu=mp+r$ and $a^\nu=np+r$, it will be that $a^\mu-a^\nu=(m-n)p$,
and thus that the difference $a^\mu-a^\nu$ of these terms will be divisible by
$p$. Therefore the difference between two terms of the given geometric
progression will be divisible by the number $p$ in innumerably many ways. 

\begin{center}
{\Large Corollary 4.}
\end{center}

6. If the power $a^\mu$ were to give the residue $r$ and the power
$a^\nu$ the residue $s$, and it were $r+s=p$, we say in this case of the
residues $r$ and $s$ for each to be the complement of the other, in
which case the sum of the powers $a^\mu+a^\nu$ will be divisible by $p$.
For if it were $a^\mu=mp+r$ and $a^\nu=np+s$, it will then be $a^\mu+a^\nu=
(m+n)p+r+s=(m+n+1)p$, and thus it has $p$ as a factor.

\begin{center}
{\Large Theorem 2.}
\end{center}

7. If the power $a^\mu$ divided by $p$ has the residue $r$, and the power
$a^\nu$ the residue $s$, the power $a^{\mu+\nu}$ will have the residue $rs$.

\begin{center}
{\Large Demonstration.}
\end{center}

Were it $a^\mu=mp+r$ and $a^\nu=np+s$, it will be
$a^{\mu+\nu}=mnpp+mps+npr+rs$; and thus if
$a^{\mu+\nu}$ is divided by $p$, the residue will be $rs$; and if it were
larger than $p$, by subtracting $p$ as many times as possible, it
can be reduced to a residue less than $p$. Q.E.D.

\begin{center}
{\Large Corollary 1.}
\end{center}

8. Hence the residue of the base $a$ divided by $p$ can be expressed
by $a$ (for if $a<p$, $a$ will be a residue properly called, and if
however $a>p$,
nevertheless the residue can be expressed by $a$, because
$a-p$
or $a-np$ are likewise treated), if the residue of the power
$a^\mu$ divided by $p$ were $r$, the residue of the power $a^{\mu+1}$ will
be $ar$, and in the same way the residue of the power $a^{\mu+2}$ will be
$a^2r$,
$a^{\mu+3}$ will be $a^3r$, etc.

\begin{center}
{\Large Corollary 2.}
\end{center}

9. From this it also follows that if the residue of the power $a^\mu$ divided by
$p$ is $=r$, for the residue of the power $a^{2\mu}$ to be
$=rr$, the residue of the power $a^{3\mu}$ to be $=r^3$, etc.
Thus if the residue of the power $a^\mu$ is $=1$, 1 will be the identical
residue of all the powers $a^{2\mu},a^{3\mu},a^{4\mu},a^{5\mu}$, etc.

\begin{center}
{\Large Corollary 3.}
\end{center}

10. But if the residue of the power $a^\mu$ divided by $p$ is $=p-1$,
which we have seen is able to be expressed by $-1$, then the residue
of the power $a^{2\mu}$ will be $=+1$, the residue of
the power $a^{3\mu}$ will be $=-1$, and in turn the residue of the power $a^{4\mu}$
$=+1$. And thus in general the residue of the power $a^{n\mu}$
will be either $+1$, if $n$ is an even number, or $-1$, if $n$ is an odd number.

\begin{center}
{\Large Scholion.}
\end{center}

11. Hence the way 
is inferred for the residues to be easily found which remain from the
division
of any power by a given number.
For say we want to investigate the residue which arises from the division
of the power $7^{160}$ by 641.

\begin{tabular}{p{1in}|p{1in}|p{4in}}
powers&residues&since of course the first power of 7 leaves 7,\\
$7^1$&7&the powers $7^2,7^3,7^4$ leave 49,343 and 478, or $-163$;\\
$7^2$&49&the square of this, $7^8$, leaves $163^2$ or $288$, and the square\\
$7^3$&343&of this, $7^{16}$, leaves $288^2$ or 255. In the same way\\
$7^4$&478&the power $7^{32}$ leaves $255^2$ or 284, and\\
$7^8$&288&the residue of the power $7^{64}$ will be $-110$, and from\\
$7^{16}$&255&$7^{128}$ arises $110^2$ or $-79$, which residue is\\
$7^{32}$&284&multiplied by 288, which gives the residue\\
$7^{64}$&$-110$&of the power $7^{128+32}=7^{160}$,\\
$7^{128}$&$-79$&which will be 640 or $-1$\\
$7^{160}$&$-1$&
\end{tabular}
We know therefore that if the power $7^{160}$ is divided by 641 then the residue
will be 640, or $-1$, from which we conclude for the residue of the power
$7^{320}$ to be $+1$. Therefore in general the residue of the power
$7^{160n}$ divided by 641 will be either $+1$, if $n$ is an even number,
or $-1$, if $n$ is an odd number.

\begin{center}
{\Large Theorem 3.}
\end{center}

12. If the number $a$ is prime to $p$, and the geometric progression
$1$; $a$; $a^2$; $a^3$; $a^4$; $a^5$; $a^6$; $a^7$; etc were formed, innumerably many terms occur
in it which leave the residue 1 when divided by $p$, and the exponents
of these terms constitute an arithmetic progression.

\begin{center}
{\Large Demonstration.}
\end{center}

Since the number of terms is infinite, but on the other hand no more than $p-1$
different residues can arise, it is necessary that several, and indeed
infinitely
many, terms produce the same residue $r$. Were $a^\mu$ and $a^\nu$ two such
terms leaving the same residue $r$, $a^\mu-a^\nu$ will be divisible by $p$.
As $a^\mu-a^\nu=a^\nu(a^{\mu-\nu}-1)$, 
and since the product is divisible by $p$ yet the one factor $a^\nu$ is prime to
$p$, it is necessary that the other factor $a^{\mu-\nu}-1$ is divisible by
$p$; from this, the power $a^{\mu-\nu}$ when divided by $p$ will have a
residue $=1$. If it were $\mu-\nu=\lambda$, so that the residue of
the power $a^\lambda$ were $=1$, it will then likewise occur
that all the powers $a^{2\lambda},a^{3\lambda},a^{4\lambda},a^{5\lambda}$, etc.
will have the same residue $=1$. And thus unity will be the residue of
all of these powers:
\[
1; a^\lambda; a^{2\lambda}; a^{3\lambda}; a^{4\lambda}; a^{5\lambda};
a^{6\lambda}; \textrm{etc.}
\]
whose exponents advance as an arithmetic progression.

\begin{center}
{\Large Corollary 1.}
\end{center}

13. Therefore if a single power $a^\lambda$ is found which when divided
by $p$ produces a residue $=1$, from it infinitely many
other powers can be exhibited which leave unity when divided by $p$. And indeed
also,
the least power of this type is $a^0=1$.

\begin{center}
{\Large Corollary 2.}
\end{center}

14. But even if no power of $a$ aside from unity appears, that is which when
divided by $p$ leaves unity as its residue, we still find for infinitely
many powers of this type to be given.

\begin{center}
{\Large Corollary 3.}
\end{center}

15. From the demonstration it is thus apparent how to give a power
$a^\lambda$ which presents a residue $=1$, whose exponent $\lambda$
is less than $p$. For if the geometric progression is continued up to the
term $a^{p-1}$, since the number of terms is $p$, it is necessary that
 at least two terms, say $a^\mu$ and $a^\nu$, have the same residue;
 from this the power $a^{\mu-\nu}$ will have a residue $=1$,
 and because $\mu<p$ and $\nu<p$, it is certain to be $\mu-\nu<p$.

\begin{center}
{\Large Theorem 4.}
\end{center}

16. If the power $a^\mu$ divided by $p$ leaves a residue $=r$, and 
the residue of a higher power $a^{\mu+\nu}$ is $=rs$, the residue
of the power $a^\nu$, with which this was increased, will be $=s$.

\begin{center}
{\Large Demonstration.}
\end{center}

As the power $a^\nu$ gives a residue, which may be set $=t$,
and since the residue of the power $a^\mu$ is $=r$, the residue of
the power $a^{\mu+\nu}$ will be $=rt$, which should be 
equal to $rs$ itself. Therefore it is $rt=rs+np$, supposing that we  set the
residues $r,t$ to be less than the divisor $p$. Therefore it shall be $t=s+
\frac{np}{r}$: but since $a$ and $p$ are mutually prime, all the residues which
arise from powers of $a$ divided by $p$ will likewise be prime to $p$,
except in the case they are $=1$. Therefore for $\frac{np}{r}$ to be an
integral
number, it is necessary that $\frac{n}{r}$ be an integral number;
putting this $=m$, it will therefore be $t=s+mp$ and thus
$t=s$. Whence if the residue of the
power $a^\mu$ is $=r$ and the residue of the power
$a^{\mu+\nu}$ is $=rs$, then it follows that the residue
of the power $a^\nu$ will be $=s$.

\begin{center}
{\Large Corollary 1.}
\end{center}

17. Therefore if $s=1$, that is if the two powers $a^\mu$ and $a^{\mu+\nu}$ have the
same residue $r$, it will follow that if the larger is divided by the smaller,
the quotient $a^\nu$ will possess a residue $=1$, the demonstration
of which rests on the preceding theorem.

\begin{center}
{\Large Corollary 2.}
\end{center}

18. If $r=1$ and $s=1$, that is if the two powers $a^\mu$ and $a^{\mu+\nu}$ have
the same
residue $=1$, then also the power $a^\nu$, whose exponent
is the difference of the previous exponents, will likewise have a residue $=1$. 

\begin{center}
{\Large Scholion.}
\end{center}

19. The theorem can also be disclosed in this way. As $a^\mu$ leaves
the residue $r$ when divided by $p$, and in the same way $a^{\mu+\nu}=np+rs$;
then it will be $a^{\mu+\nu}-a^\mu s=np-mps=(n-ms)p$; and thus the number
$a^{\mu+\nu}-a^\mu s=a^\mu(a^\nu-s)$ will be divisible by $p$:
but since the first factor $a^\mu$ is not divisible by $p$, therefore
the other factor $a^\nu-s$ will be divisible by $p$, and it follows that the
power $a^\nu$ gives a residue $=s$ when divided by $p$.

\begin{center}
{\Large Theorem 5.}
\end{center}

20. If after unity $a^\lambda$ is the smallest power which leaves
unity when divided by $p$, then no other powers may leave the same residue
$=1$, unless they occur in this geometric progression.
\[
1; a^\lambda; a^{2\lambda}; a^{3\lambda}; a^{4\lambda}; a^{5\lambda};
\textrm{etc.} 
\]

\begin{center}
{\Large Demonstration.}
\end{center}

For we may put some other power $a^\mu$ to give a residue $=1$
when divided by $p$, and with $\mu>\lambda$, neither however is this
equal to any multiple of $\lambda$. Thus this exponent $\mu$ can be expressed
as $\mu=n\lambda+\delta$, where $\delta<\lambda$, and neither will it
be $\delta=0$. Therefore since the power $a^{n\lambda}$, and to
the same extent $a^\mu=
a^{n\lambda+\delta}$, will leave unity when divided by $p$, by \S 18 the power
$a^\delta$ will have unity as its residue, and it would therefore
be that $a^\lambda$ was not the smallest power of this type, contrary to the
hypothesis. Whence if $a^\lambda$ is the smallest power giving a residue
$=1$, then no other powers will have this same property except those whose
exponents are multiples of $\lambda$.

\begin{center}
{\Large Corollary 1.}
\end{center}

21. Therefore if the second term of the geometric progression $1,a,a^2,a^3,
a^4$, etc. leaves 1 when divided by $p$, then all the terms will yield the
same residue $=1$: therefore no residue aside from 1 will occur
in any of these numbers.

\begin{center}
{\Large Corollary 2.}
\end{center}

22. If the residue of the third term $a^2$ is $=1$,
that is if $a^2=np+1$, then the alternate terms $1,a^2,a^4,a^4$, etc.,
whose exponents are even, will all have the same residue $=1$,
and indeed all the other terms, unless $a^1$ also has a residue $=1$,
will all produce different residues.

\begin{center}
{\Large Corollary 3.}
\end{center}

23. It can happen, therefore, that many 
fewer numbers may occur among the residues than the number $p - 1$ 
contains 
unities: on the other hand, more than $p - 1$ 
distinct numbers cannot occur.

\begin{center}
{\Large Theorem 6.}
\end{center}

24. If the power $a^{2n}$, whose exponent is an even number, divided by a
prime number $p$ leaves a residue $=1$, then the power $a^n$ divided by
the same number $p$ will give a residue either $=+1$ or $=-1$.

\begin{center}
{\Large Demonstration.}
\end{center}

Indeed, we put $r$ to be a residue which remains from the division of the
power $a^n$ by the prime number $p$, and the residue of the power $a^{2n}$ will
be $=rr$, which by the hypothesis is $=1$. From this it will
be that $rr=1+mp$, and so $rr-1=mp$; from which since $rr-1=(r+1)(r-1)$ must
be divisible by $p$, either the factor $r+1$ or the factor $r-1$ shall be
divisible by $p$. In the first case it will be $r+1=\alpha p$, and $r=\alpha
p-1$, and hence $r=-1$. In the latter case it will be $r-1=ap$ and $r=ap+1$
and hence $r=+1$. Therefore if the power $a^{2n}$ yields a residue
$=+1$, the power $a^n$ will have a residue $=+1$
or $=-1$, if indeed $p$ were a prime number.

\begin{center}
{\Large Corollary 1.}
\end{center}

25. Consequently if $a^{2n}$ were the minimum power which leaves a residue
$=+1$ when divided by the prime number $p$, then the power $a^n$ will give
a residue $=-1$. Therefore if the exponent $\lambda$ of the minimum
power $a^\lambda$ yielding a residue $=1$ is even, then among the residues
of the terms of the geometric progression $1,a,a^2,a^3,a^4$, etc. the number
$-1$ will also appear.

\begin{center}
{\Large Corollary 2.}
\end{center}

26. But if on the other hand the exponent $\lambda$ of the minimum power
$a^\lambda$ yielding the residue 1 were an odd number, then no power
whatsoever will leave a residue $=-1$. For if $a^\mu$ were such a power,
which would give a residue $=-1$, then the power
$a^{2\mu}$ would give a residue $=+1$, whence $2\mu=n\lambda$, and
since $\lambda$ is an odd number it would be $2\mu=2m\lambda$, and so
$\mu=m\lambda$. But then the power $a^{m\lambda}$ leaves a residue
$=+1$,
and therefore the residue $-1$ does not occur anywhere.

\begin{center}
{\Large Theorem 7.}
\end{center}

27. If $a^\lambda$ is the minimum power of $a$ which yields a residue $=1$
when divided by the number $p$, then all the residues which result from
the terms of the geometric progression $1,a,a^2,a^3,\ldots,a^{\lambda-1}$,
continued onto the power $a^\lambda$, will be mutually distinct.

\begin{center}
{\Large Demonstration.}
\end{center}

For if, say, $a^\mu$ and $a^\nu$ are two powers whose exponents $\mu$ and $\nu$
are less than $\lambda$ and which give the same
residue, the the difference $a^\mu-a^\nu$ of them would be divisible by
$p$, and thus the power $a^{\mu-\nu}$ divided by $p$ would leave a residue 
$=+1$, for which it would be $\mu-\nu<\lambda$, contrary to the
hypothesis;
from this it is clear for all the powers whose exponents are less than $\lambda$
to yield different residues.

\begin{center}
{\Large Theorem 8.}
\end{center}

28. If $a^\lambda$ were a certain power of $a$, which divided by the number
$p$ produces a residue $=1$, and also the geometric progression
were divided into sections according to the powers
$a^\lambda,a^{2\lambda},a^{3\lambda},a^{4\lambda}$, etc. in this way:
\[
1,a,a^2,\ldots,a^{\lambda-1}|a^\lambda,\ldots,a^{2\lambda-1}|a^{2\lambda},\ldots,
a^{3\lambda-1}|a^{3\lambda},\ldots,a^{4\lambda-1}|\, \textrm{etc.}
\]
such that each section contains $\lambda$ terms, then  
in each section the same residues will appear,
and also the same order will be repeated.

\begin{center}
{\Large Demonstration.}
\end{center}

Indeed, the first terms $1,a^\lambda,a^{2\lambda},a^{3\lambda}$, etc. of
all the sections produce the same residue $=1$.
Then the second terms $a,a^{\lambda+1},a^{2\lambda+1},a^{3\lambda+1}$, etc.
of all the sections will likewise give the same residue;
for if $r$ were a residue arising from the term $a^1$, since
$a^{\lambda+1}=a^\lambda \cdot a^1$, the residue arising from this term will be equal
to $1\cdot r=r$; and in a similar way it is apparent for the residues
of all the terms $a^{2\lambda+1},a^{3\lambda+1}$, etc. to be $=r$.
And if in general $a^\mu$ were a certain term of the first section,
the residue arising from which were $=r$, then too the residue of
the term $a^{n\lambda+\mu}$ will be $=r$, since the residue
of the term $a^{n\lambda}$ is $=1$: hence the analogous terms
$a^{\lambda+\mu},a^{2\lambda+\mu},a^{3\lambda+\mu}$, etc. of all
the sections will have the same residue.

\begin{center}
{\Large Corollary 1.}
\end{center}

29. And if so many many of the terms contained in the first section were known,
then too the residues of all the terms which constitute the remaining sections
will be known.

\begin{center}
{\Large Corollary 2.}
\end{center}

30. For if a term $a^x$ is given, whose exponent $x$ is an arbitrarily
large number, the residue of it can be easily found. For this exponent
$x$ can be reduced to the form $n\lambda+\mu$, such that $\mu<\lambda$,
and thus the residue of the term $a^x$ will be the same as the that of
the term $a^\mu$.

\begin{center}
{\Large Corollary 3.}
\end{center}

31. Moreover for this number $\mu$ which is chosen as less than $\lambda$,
if the number $x$ is divided by $p$, then for the residue   which remains from this
division which is sought will be this very number $\mu$.

\begin{center}
{\Large Corollary 4.}
\end{center}

32. Moreover powers $a^\lambda$ may always be given which when divided by
$p$ leave unity, whose exponent $\lambda$ is less than the given
number $p$, and thus for finding the residues of all the terms of the geometric
progression, it is not necessary to continue beyond the term $a^p$.

\begin{center}
{\Large Corollary 5.}
\end{center}

33. If moreover the power $a^\mu$ is the minimum of those which leave unity
when divided by the number $p$, then as all the terms less than $a^\lambda$
yield
different residues, in all the residues, no more and no less different
numbers occur than $\lambda$. Therefore if $\lambda$ were less than $p-1$,
not all numbers occur in the residues: rather certain numbers will
clearly never be able to remain in the division of the terms of the geometric
progression $1,a,a^2,a^3$, etc.

\begin{center}
{\Large Corollary 6.}
\end{center}

34. Therefore if the diversity of the residues is considered, it can happen
that only a single residue occurs among all the powers of $a$, or
exactly two different residues, or three etc. are produced, though
no more than $p-1$ can occur. But however many residues should be produced,
unity will always appear among them.

\begin{center}
{\Large Theorem 9.}
\end{center}

35. If $p$ were a prime number and $a$ were prime to $p$, and also
all the numbers less than $p$ itself appear among the residues which 
arise from the division of all the powers of $a$ by the number $p$, then
$a^{p-1}$ will be the minimum power which leaves unity when divided by $p$.

\begin{center}
{\Large Demonstration.}
\end{center}

Were $a^\lambda$ the minimum power which leaves unity when divided by $p$,
it is clear from the preceding that $\lambda < p$ (15). Now since the number
of all the residues of divisors is $=\lambda$, and all the numbers
less than $p$ is $=p-1$, it is clear that if it were $\lambda<p-1$, then not all
the numbers less than $p$ could occur in the residues; therefore it will
not be $\lambda<p-1$, and neither indeed is it $\lambda>p-1$, since in general
it
is $\lambda<p$. From this it follows that $\lambda=p-1$. 
Whenceforth if all the numbers less than $p$ occur in the residues, the
power $a^{p-1}$ will be the minimum which leaves unity when divided by $p$.

\begin{center}
{\Large Scholion.}
\end{center}

36. The nature of this theorem requires that $p$ be a prime number; for
unless this is the case, it is not possible to make it that all the numbers
less than $p$ occur in the residues. For it may be clearly seen, which
is to be carefully considered,
that if $p$ is a composite number, to which however $a$ is prime, for no aliquot
part of $p$ to occur in the residues: for if some power $a^\mu$ should give a
residue $r$ which was an aliquot part of $p$, as $a^\mu=mp+r$, 
this power $a^\mu$ would have the divisor $r$, and thus neither it
nor its root $a$ would be prime to $p$, which is against the hypothesis.

\begin{center}
{\Large Theorem 10.}
\end{center}

37. If the number of different residues which arise from the division
of the powers $1,a,a^2,a^3,a^4,a^5$, etc. by a prime number $p$ is less
than $p-1$, then there are at least as many numbers which are not residues
as which are residues.

\begin{center}
{\Large Demonstration.}
\end{center}
Were $a^\lambda$ the minimum power which leaves unity when divided by $p$,
and also $\lambda<p-1$, the number of all the residues of the divisor
will be $=\lambda$, and thus less than $p-1$. Therefore as the number of all
the numbers less than $p$ itself are $=p-1$, it is clear in the given
case for numbers to be given which do not occur as residues. I say moreover
for the number of numbers of this type at a minimum to be $=\lambda$. So that
this may be exhibited, we shall relate the residues by the terms from which
they arise,
and they will be 
\[
\textrm{these residues} \quad 1,a,a^2,a^3,a^4,\ldots,a^{\lambda-1}
\]
the number of which is $=\lambda$, and these residues, if reduced to the usual
form, will all be less than $p$ and mutually distinct. Therefore
as it was $\lambda<p-1$ by hypothesis, a number is certainly given which
does not appear in these residues. Let $k$ be such a number; I say now
that if $k$ is not a residue, then neither $ak$, nor $a^2k$, nor $a^3k$, etc.,
nor $a^{\lambda-1}k$ will occur in the residues. For
let $a^\mu k$ be a residue arising from the power $a^\alpha$, and it shall be
$a^\alpha=np+a^\mu k$, that is $a^\alpha-a^\mu k=np$, and thus $a^\alpha-a^\mu
k=a^\mu (a^{\alpha-\mu}-k)$ is divisible by $p$. But since $a^\mu$ is not
divisible by $p$, $a^{\alpha-\mu}$ will therefore be divisible by $p$, that is
the power $a^{\alpha-\mu}$ when divided by $p$ will leave the residue $k$,
which opposes the hypothesis. Whence it is clear that all these numbers:
$k,ak,a^2k,a^3k,\ldots,a^{\lambda-1}k$, that is all the numbers derived from
it, to not be residues. And these numbers, whose multitude is $=\lambda$,
are all mutually distinct; for if two, say $a^\mu k$ and $a^\nu k$, convene,
and reduce to the same residue $r$, it would
be $a^\mu k=mp+r$ and $a^\nu k=np+r$, and thus
$a^\mu k-a^\nu k=(m-n)p$, that is $(a^\mu-a^\nu)k=(m-n)p$ would be divisible
by $p$. Neither indeed is $k$ divisible by $p$, as we took $p$ as a prime number
and $k<p$; were $a^\mu-a^\nu$ divisible by $p$, or should $a^{\mu-\nu}$ leave unity
when divided by $p$, since however $\mu<\lambda-1$ and $\nu<\lambda-1$, it would
be $\mu-\nu<\lambda$, which is absurd. Therefore all these numbers
$k,ak,a^2k,a^3k,\ldots,a^{\lambda-1}k$, if reduced, will be mutually
distinct, and the multitude of them is $=\lambda$. 
Therefore at a minimum $\lambda$ numbers are given which will not occur in the
residues, supposing that it were $\lambda<p-1$.

\begin{center}
{\Large Corollary 1.}
\end{center}

38. Consequently were there $\lambda$ different numbers which are residues,
and just as many different numbers which are not residues, and all were
smaller than $p$, the number of them joined together cannot be larger than
$2\lambda$, or $p-1$: for there are no more numbers less than $p$ given
than $p-1$.

\begin{center}
{\Large Corollary 2.}
\end{center}

39. Therefore if $a^\lambda$ were the minimum power which leaves unity
when divided by the prime number $p$, it would be $\lambda<p-1$, and then
it will certainly not be that $\lambda>\frac{p-1}{2}$: it will therefore
be
either $\lambda=\frac{p-1}{2}$ or $\lambda<\frac{p-1}{2}$.

\begin{center}
{\Large Corollary 3.}
\end{center}

40. From before we see for the exponent $\lambda$ of the minimum power
to necessarily be less than $p$; it will therefore be either
$\lambda=p-1$ or $\lambda<p-1$; in the case that
$\lambda<p-1$, we know at once for it now to be either $\lambda=\frac{p-1}{2}$
or $\lambda<\frac{p-1}{2}$. And thus no number is contained beyond the
limits $p-1$ and $\frac{p-1}{2}$ which could ever be a value of $\lambda$.

\begin{center}
{\Large Theorem 11.}
\end{center}

41. If $p$ were a prime number, and also $a^\lambda$ were the minimum power
of $a$ which leaves unity when divided by $p$, and it were
$\lambda < \frac{p-1}{2}$; then it cannot happen that the exponent $\lambda$
is larger than $\frac{p-1}{3}$: and it will therefore be either $\lambda=
\frac{p-1}{3}$ or $\lambda<\frac{p-1}{3}$.

\begin{center}
{\Large Demonstration.}
\end{center}

Were $a^\lambda$ the minimum power which leaves unity when divided by a prime
number $p$, then no more than $\lambda$ different numbers occur in the residues,
which are left from this terms
\[
1; a; a^2; a^3; a^4; \ldots; a^{\lambda-1}
\]
if each were divided by $p$; whence since $\lambda<p-1$, $p-1-\lambda$ numbers
will occur which are not residues, and if one of these is $=r$, we have seen
that all these numbers
\[
r; ar; a^2r; a^3r; a^4r; \ldots; a^{\lambda-1}r,
\]
according as they are reduced to numbers less than $p$ itself when divided by $p$, 
to not be contained in the residues. Then moreover at least $\lambda$ numbers
are excluded from the residues; whence as it was $\lambda<\frac{p-1}{2}$, it will
be $\lambda<p-1-\lambda$, and thus aside from these numbers others are also
given which are not contained in the residues. Let $s$ be a number of this
type, which was neither a residue or contained in the preceding series of
non-residues; then too all these numbers
\[
s; as; a^2s; a^3s; a^4s; \ldots; a^{\lambda-1}s
\]
will not be residues: and these numbers, which we have displayed
in the preceding demonstration, will all be mutually distinct. Neither indeed
shall any of these numbers, say $a^\mu s$, be contained now in the preceding
series of non-residues, that is it is not $a^\mu s=a^\nu r$. For indeed if
it were $a^\nu r=a^\mu s$, it would become $s=a^{\nu-\mu}r$, or
$s=a^{\lambda+\nu-\mu}r$, supposing that it were $\mu>\nu$, from which
now $s$ shall be contained in the prior series contrary to the hypothesis. 
Wherefore if $\lambda<\frac{p-1}{2}$, 
at least $\lambda$ numbers are given which are not residues,
and since we have $\lambda$ residues and $2\lambda$ non-residues, 
and all these numbers were less than $p$ itself,
it cannot happen that the sum $3\lambda$ of them would be greater than $p-1$,
that
is it will not be that $\lambda>\frac{p-1}{3}$.
Therefore it will be either $\lambda=\frac{p-1}{3}$ or $\lambda<\frac{p-1}{3}$,
if indeed $\lambda<\frac{p-1}{2}$ and $p$ were a prime number.

\begin{center}
{\Large Corollary 1.}
\end{center}

42. Therefore if it is not $\lambda<\frac{p-1}{3}$, then it will certainly
be $\lambda=\frac{p-1}{3}$, just as when $\lambda<\frac{p-1}{2}$. And 
by removing this condition, if we know that it is not $\lambda<\frac{p-1}{3}$,
then it necessarily follows that either $\lambda=\frac{p-1}{3}$,
$\lambda=\frac{p-1}{2}$, or $\lambda=p-1$.

\begin{center}
{\Large Corollary 2.}
\end{center}

43. Moreover if it were either $\lambda=\frac{p-1}{3}$ or
$\lambda=\frac{p-1}{2}$,
the power $a^{p-1}$ divided by $p$ will leave unity. For if $a^\lambda$ leaves
unity, then also $a^{2\lambda}$ and $a^{3\lambda}$ will give unity as residues.

\begin{center}
{\Large Theorem 12.}
\end{center}

44. If $a^\lambda$ were the minimum power of $a$ which leaves unity when divided
by a prime number $p$, and it were $\lambda<\frac{p-1}{3}$, then it will
certainly not be $a>\frac{p-1}{4}$, and therefore it will be either
$\lambda=\frac{p-1}{4}$ or $\lambda<\frac{p-1}{4}$.

\begin{center}
{\Large Demonstration.}
\end{center}

Because the number of all the different residues which 
come forth from the division of all the powers of $a$ by the prime number $p$
is $=\lambda$,  
and in fact arise from these terms: $1; a; a^2; a^3; a^4;\ldots; a^{\lambda-1}$:
on account of $\lambda<\frac{p-1}{3}$, immediately twice as many
numbers will occur which are not residues, which arise from these two
progressions
\begin{eqnarray*}
&&r;ar;a^2r; a^3r; a^4r; \ldots; a^{\lambda-1}r\\
\textrm{and}&&s;as;a^2s;a^3s;a^4s; \ldots; a^{\lambda-1}s
\end{eqnarray*}
the number of these numbers, both residues and non-residues, is $=3\lambda$,
and thus less than $p-1$, therefore numbers will still be left over which will
not be residues.
Were $t$ such a number, we display as before 
that too all these numbers
\[
t; at; a^2t; a^3t; a^4t; \ldots; a^{\lambda-1}t
\]
will not be residues, the number of which is $=\lambda$. And these numbers
are not only all mutually distinct, since $p$ was a prime number, but moreover
they are different from the preceding, and thus the multitude of all
these numbers, either residues or non-residues, is $=4\lambda$, and also
all these numbers shall be less than $p$; it is impossible that
it were $4\lambda>p-1$; therefore it will be either
$\lambda=\frac{p-1}{4}$ or $\lambda<\frac{p-1}{4}$, if indeed, as we have
assumed,
$\lambda<\frac{p-1}{3}$ and $p$ were a prime number.

\begin{center}
{\Large Corollary 1.}
\end{center}

45. In a similar way it can be demonstrated, that if $\lambda<\frac{p-1}{4}$
then it is impossible that $\lambda>\frac{p-1}{5}$, and to thus be either
$\lambda=\frac{p-1}{5}$ or $\lambda<\frac{p-1}{5}$.

\begin{center}
{\Large Corollary 2.}
\end{center}

46. In general if it is known that $\lambda<\frac{p-1}{n}$, it may be
demonstrated
in the same way that it cannot be $\lambda>\frac{p-1}{n+1}$, and therefore
will be either $\lambda=\frac{p-1}{n+1}$ or $\lambda<\frac{p-1}{n+1}$.

\begin{center}
{\Large Corollary 3.}
\end{center}

47. From this it is clear for the number of all those numbers which cannot
be residues to be either $=0$, or $=\lambda$, or $=2\lambda$, or another
such multiple of $\lambda$: for if they were more numbers of this type
than $n\lambda$, then because from each $\lambda$ others follows, the number
of all of these would be $=(n+1)\lambda$; and if these are not
the only numbers contained in the non-residues, then again $\lambda$
others follow simultaneously.

\begin{center}
{\Large Theorem 13.}
\end{center}

48. If $p$ were a prime number, and $a^\lambda$ the minimum power of $a$
which leaves unity when divided by $p$, the exponent $\lambda$ will be a divisor
of the number $p-1$.

\begin{center}
{\Large Demonstration.}
\end{center}

Therefore the number of all the residues of divisors is $=\lambda$,
from which the number of the remaining numbers less than $p$ which are
not residues will be $=p-1-\lambda$; but this number (47) is a multiple
of $\lambda$, say put as $n\lambda$, thus it would be that
$p-1-\lambda=n\lambda$, from which it would be $\lambda=\frac{p-1}{n+1}$. 
It is therefore evident for the exponent $\lambda$ to be a divisor of the number
$p-1$, from which if it were not $\lambda=p-1$, it will be certain
for the exponent $\lambda$ to be equal to some aliquot part
of the number $p-1$.

\begin{center}
{\Large Theorem 14.}
\end{center}

49. If $p$ were a prime number, and $a$ prime to $p$, then the power
$a^{p-1}$ will leave unity when divided by $p$.

\begin{center}
{\Large Demonstration.}
\end{center}

Were $a^\lambda$ the minimum power of $a$ which leaves unity when divided by
$p$,
as we have seen $\lambda<p$, and indeed we have demonstrated above
for it to be either that $\lambda=p-1$ or for $\lambda$ to be an aliquot
part of the number $p-1$. In the first case
the proposition is obvious, as the power $a^{p-1}$ leaves unity when 
divided by $p$. In the latter case, in which $\lambda$ is an aliquot
part of the number $p-1$, it will be $p-1=n\lambda$, but since the power
$a^\lambda$ leaves unity when divided by $p$, so too all the powers
$a^{2\lambda},a^{3\lambda},a^{4\lambda}$, etc. and therefore $a^{n\lambda}$,
that is $a^{p-1}$, shall leave unity when divided by $p$. Therefore
the power $a^{p-1}$ will always leave unity when divided by $p$.

\begin{center}
{\Large Corollary 1.}
\end{center}

50. Since the power $a^{p-1}$ leaves unity when it is divided by the prime
number $p$, the formula $a^{p-1}-1$ will be divisible by the prime number $p$,
supposing that $a$ were prime to $p$, that is that $a$ were not divisible by
$p$.

\begin{center}
{\Large Corollary 2.}
\end{center}

51. Therefore if $p$ were a prime number, all the powers of the exponent
$p-1$, such as $n^{p-1}$, when divided by $p$ will leave either unity or nothing.
The first happens if $n$ is a number prime to $p$, the latter
indeed if $n$ is a number divisible by $p$.

\begin{center}
{\Large Corollary 3.}
\end{center}

52. If $p$ were a prime number, and $a$ and $b$ were numbers prime to $p$,
the difference of the powers $a^{p-1}-b^{p-1}$ will be divisible by the number
$p$. For since both $a^{p-1}-1$ and $b^{p-1}-1$ are divisible by $p$, then too
the differences of these formulas, that is $a^{p-1}-b^{p-1}$, will be divisible
by $p$.

\begin{center}
{\Large Scholion.}
\end{center}

53. Behold therefore a new demonstration of the extraordinary theorem,
found before by Fermat, which differs greatly from that which I gave
in the Comment. Acad. Petropol. Volume VIII. For there
I called upon the expansion of the binomial $(a+b)^n$ into a series
by means of the method of Newton, which reasoning seems quite remote
from the proposition; here indeed I have demonstrated the same theorem
from the properties of powers alone, by which this demonstration seems
much more natural. In addition, other important properties
about the residues of powers when divided by prime numbers may appear
to us.
For it is clear that if $p$ were a prime number, then not only will the formula
$a^{p-1}-1$ be divisible by $p$, but it will also sometimes happen that the
simpler
formula $a^\lambda-1$ will be divisible by $p$, and then for the exponent
$\lambda$ to be an aliquot part of the exponent $p-1$.

\begin{center}
{\Large Theorem 15.}
\end{center}

54. If $q$ were a prime number, and also the power $a^q$ left unity when divided
by a prime number $p$, then $a^q$ will be the minimum power of $a$ which leaves
unity when divided by $p$, except for the singular case in which the number
$a$ leaves unity when divided by $p$.

\begin{center}
{\Large Demonstration.}
\end{center}

For were $a^\lambda$ the minimum power of $a$ which left unity when divided by
the prime number $p$,
then no other powers will be gifted with this property, except
for $a^{2\lambda},a^{3\lambda},a^{4\lambda}$,  etc. Truly none
of these will be able to be equal to the power $a^q$, unless it were
$\lambda=1$, since $q$ was a prime number, and thus it is necessary that it
be $q=\lambda$, and thus that $a^q$ is the minimum power which leaves
unity when divided by $p$. However the case in which $\lambda=1$ is excepted,
that is
in which the number $a$ leaves unity when divided by $p$: for in this case
each power $a^n$, either when the exponent $n$ was a prime number or a composite
number, shall leave unity when division by $p$ is carried out.

\begin{center}
{\Large Corollary 1.}
\end{center}

55. Therefore if the power $a^q$, whose exponent is a prime number, leaves unity
when divided by a prime number $p$, then $q$ will be an aliquot part
of the number $p-1$, in which case the formula $a^q-1$ will be divisible by
the prime number $p$.

\begin{center}
{\Large Corollary 2.}
\end{center}

56. Were $q$ an aliquot part of the number $p-1$, it will be
$p-1=nq$, and $p=nq+1$. And if therefore the formula
$a^q-1$, in which $q$ is a prime number, were divisible by this prime
number $p$, it will always have to have a divisor of such form $p=nq+1$,
unless $p=a-1$:
for $a-1$ is always a divisor of the formula $a^q-1$.

\begin{center}
{\Large Corollary 3.}
\end{center}

57. Therefore the formula $a^q-1$, with $q$ a prime number, will
not admit prime divisors other than $p-1$ unless they are contained in the
form $nq+1$; and as $q$ is a prime number, and therefore odd, unless
$q=2$, except for which even numbers can be taken for $n$,
and therefore all the divisors,
if it has any, will be contained
in the form $2nq+1$.

\begin{center}
{\Large Corollary 4.}
\end{center}

58. Therefore since a divisor of the formula $a^q-1$ is
\[
a^{q-1}+a^{q-2}+a^{q-3}+a^{q-4}+\cdots+a^2+a+1
\]
this form will be contained in $2nq+1$, and therefore
this expression: $a^{q-1}+a^{q-2}+a^{q-3}+\cdots+a^2+a$ will be divisible
by the prime number $p$, whatever number $a$ is, though if $a=q$ or $a=mq$ this
is manifest by itself.

\begin{center}
{\Large Scholion 1.}
\end{center}

59. This is moreover clear, if $a$ were neither $q$ nor $mq$; then for the 
formula which has been found to turn into 
\[
a(a^{q-2}+a^{q-3}+a^{q-4}+\cdots+a+1)
\]
of which the second factor, which transforms into $\frac{a^{q-1}-1}{a-1}$,
is divisible by $q$: because indeed by this it is evident; 
for since $q$ was a prime number, the formula $a^{q-1}-1$ is divisible by it;
and when this is further divided by $a-1$, it will remain
divisible by $p$, unless $a-1$
has $q$ as a divisor, which case is excepted now as before. For it is
to be noted for the form $a^{q-1}+a^{q-2}+a^{q-3}+\cdots+a^2+a+1$ to be
contained in the form $2nq+1$, insofar as this is either a prime number, or
composed
from prime numbers of the form $2nq+1$. But if these formula
should have now a factor $a-1$, by which the form $a^q-1$ is divisible,
then the form $2nq+1$ shall not agree with it. But if $a-1=mq$, that is
$a=mq+1$, then this
form will be divisible by $q$, because the number of terms is $=q$,
and therefore this will not be contained in the form $2nq+1$.

\begin{center}
{\Large Scholion 2.}
\end{center}

60. It is also most interesting to know the divisors of the formula
$a^q-1$, when $q$ is a prime number, because usually these, except the divisor
$a-1$ which reveals itself at once, are most difficult to investigate,
and indeed it is possible to happen that often such a formula,
after it is divided by $a-1$, shall be a prime number. And if $q$  is not
a prime number, but should have divisors of itself $m,n$, then clearly
the formulas $a^m-1$ and $a^n-1$ will be divisors of the formula $a^q-1$.
Therefore in these cases the investigation of the further divisors is
reduced to the formulas $a^m-1$ and $a^n-1$, 
in which the exponents $m$ and $n$ are prime numbers. We know therefore
that if it is wanted to be attempted
to investigate the divisors of the formula $a^q-1$, an attempt should be
made with no other prime numbers, unless they were contained in the form
$2nq+1$, which operation elsewhere is most difficult,
and is not easily lightened.

\begin{center}
{\Large Theorem 16.}
\end{center}

61. If the power $a^m$, divided by the number $p$, leaves a residue $=r$,
then moreover the power $(a \pm \alpha p)^m$, divided by $p$, will leave
the same residue $r$.

\begin{center}
{\Large Demonstration.}
\end{center}

If the power $(a\pm \alpha p)^m$ is expanded, it will produce
\[
a^m \pm m\alpha a^{m-1}p\pm \frac{m(m-1)}{1\cdot 2}\alpha^2 a^{m-2}p^2 \pm
\textrm{etc.}
\]
all of whose terms, besides the first, are divisible by $p$: whence this
quantity will leave the same residue when divided by $p$ as if just the first
term $a^m$ were divided by $p$. Therefore as the power $a^m$ leaves a residue
$=r$, so too the power $(a \pm \alpha p)^m$ will leave a residue $=r$.

\begin{center}
{\Large Corollary 1.}
\end{center}

62. If $m$ were an even number, the demonstration would also be valid
for the formula $(-a+\alpha p)^m$, therefore in this case also
the formula $(\alpha p-a)^m$, when divided by $p$ will leave the same residue
$r$
that the formula $a^m$  leaves.

\begin{center}
{\Large Corollary 2.}
\end{center}

63. But if $m$ were an odd number, because the formula $-a^m$ leaves a residue
$=-r$ when divided by $p$, also the formula $(\alpha p-a)^m$ shall leave a
residue $=-r$.

\begin{center}
{\Large Theorem 17.}
\end{center}

64. If it were $a=c^n \pm \alpha p$, then the formula $a^{\frac{p-1}{n}}$
divided by the prime number $p$ shall leave unity, supposing that $n$ were a
divisor
of the number $p-1$.

\begin{center}
{\Large Demonstration.}
\end{center}

Since it is $a=c^n \pm \alpha p$, the power $a^{\frac{p-1}{n}}$,
or
$(c^n \pm \alpha p)^{\frac{p-1}{n}}$, when divided by $p$ will leave
the same residue as the power $c^{n\cdot \frac{p-1}{n}}$ or $c^{p-1}$,
and as $p$ is a prime number, the power $c^{p-1}$ will leave unity
when divided by $p$, therefore
also the power $a^{\frac{p-1}{n}}$ shall leave unity, supposing indeed
that
$a=c^n \pm \alpha p$, and neither $a$ or $c$ were divisible by $p$.

\begin{center}
{\Large Corollary 1.}
\end{center}

65. Therefore from this theorem the case may be understood in which the powers
of the numbers, whose exponents are less than $p-1$,
leave unity when divided by the prime number $p$.

\begin{center}
{\Large Corollary 2.}
\end{center}

66. Therefore if it were $a=cc+\alpha p$, with $p$ taken as a prime number,
then the power $a^{\frac{p-1}{2}}$ will leave unity when divided by $p$,
that is the formula $a^{\frac{p-1}{2}}-1$ will be divisible by $p$. For
were $p$ a prime number, unless it were $=2$, the exponent $\frac{p-1}{2}$
will always be an integral number.

\begin{center}
{\Large Corollary 3.}
\end{center}

67. If it were $a=c^3\pm \alpha p$, then the power $a^{\frac{p-1}{3}}$ will
leave
unity when divided by $p$, that is the formula $a^{\frac{p-1}{3}}-1$
will be divisible by $p$. This case will occur if the prime number $p$
is disposed such that $p-1$ were divisible by 3.

\begin{center}
{\Large Theorem 18.}
\end{center}

68. If it were $ab^n=c^n \pm \alpha p$, and $p$ were a prime number, then the
power
$a^{\frac{p-1}{n}}$ will leave unity when divided by $p$, providing that
$\frac{p-1}{n}$ is an integral number.

\begin{center}
{\Large Demonstration.}
\end{center}

The power $(c^n \pm \alpha p)^{\frac{p-1}{n}}$, or $a^{\frac{p-1}{n}}b^{p-1}$,
will leave the same residue when divided by $p$ as the power $c^{n\cdot
\frac{p-1}{n}}=c^{p-1}$, and since this power leave unity, therefore
so will the power $a^{\frac{p-1}{n}}b^{p-1}$. Moreover the factor
$b^{p-1}$ of this likewise leaves
leave unity; therefore it is necessary for the other factor
$a^{\frac{p-1}{n}}$,
if divided by $p$, to leave unity, unless either $b$ or $c$ were
divisible by $p$.

\begin{center}
{\Large Corollary 1.}
\end{center}

69. Therefore if $ab^n=c^n\pm \alpha p$, that is $ab^n-c^n$, or
$c^n-ab^n$, were divisible by the prime number $p$, then too this formula
$a^{\frac{p-1}{n}}-1$ will be divisible by $p$.

\begin{center}
{\Large Corollary 2.}
\end{center}

70. Were $p$ a prime number, where it may be put $p=mn+1$, if the formula
$ab^n-c^n$ or $c^n-ab^n$ divisible by $p$, then also the formula $a^m-1$ will be
divisible by this prime number $p$.

\begin{center}
{\Large Corollary 3.}
\end{center}

71. Therefore provided that numbers $b$ and $c$ are given of this type,
that $ab^n-c^n$ or $c^n-ab^n$ admit division by a prime number $p=mn+1$,
then it is certain for the formula $a^m-1$ to be divisible by the same
prime number $p=mn+1$.

\begin{center}
{\Large Theorem 19.}
\end{center}

72. If the formula $a^m-1$ were to be divisible by a prime number $p=mn+1$, then
there will always be given numbers $x$ and $y$ of such a type 
that $ax^n-y^n$ will be divisible by this prime number $p$.

\begin{center}
{\Large Demonstration.}
\end{center}

For since $x^{mn}$ and $y^{mn}$ leave unity when divided by $p$,
the formula $a^mx^{mn}-y^{mn}$ will always be divisible by $p$, 
provided that neither $x$ or $y$ are divisible by $p$. Since also
by factors it would be
$a^mx^{mn}-y^{mn}=(ax^n-y^n)(a^{m-1}x^{mn-n}+a^{m-2}x^{mn-2n}y^n+a^{m-3}x^{mn-3n}y^{2n}+
\cdots+y^{mn-1})$,
if this ever denies the prime factor $ax^n-y^n$ to be divisible by $p$,
it is compelled to affirm for the other
factor to always be divisible by $p$, providing that the numbers taken
for $x$ and $y$ are not divisible by $p$.  Any value whatsoever may be retained
for $x$, and for $y$ we should take the successive numbers $1,2,3,4$, onto
$p-1=mn$, not ever shall a value divisible by $p$ be obtained, and were
it for the sake of brevity
\begin{eqnarray*}
A&=&a^{m-1}x^{mn-n}+a^{m-2}x^{mn-2n}+\cdots+1\\
B&=&a^{m-1}x^{mn-n}+a^{m-2}x^{mn-2n}2^n+\cdots+2^{mn-n}\\
C&=&a^{m-1}x^{mn-n}+a^{m-2}x^{mn-2n}3^n+\cdots+3^{mn-n}\\
&\vdots&\\
N&=&a^{m-1}x^{mn-n}+a^{m-2}x^{mn-2n}(mn)^n+\cdots+(mn)^{mn-n}
\end{eqnarray*}
and also all these quantities, which constitute an algebraic progression
of order $mn-n$,  would be divisible by $p$, and then moreover
the first
differences of these, the second, third and each in this order. But the
differences of this series of order $mn-n$, which are defined only
by the $mn-n+1$ terms of this series, do not involve the term $(mn+1)^{mn-n}$,
or $p^{mn-n}$, because $p$ is not able to be a value of $y$, it follows:
\[
1\cdot 2\cdot 3\cdot 4 \cdot 5\cdots (mn-n)
\]
which plainly is not divisible
by the prime number $p=mn+1$, because its only prime divisors are less than
$mn-n$.
Therefore since the difference of this order $mn-n$ is not divisible by $p$,
it follows that not all the terms of the series $A,B,C,D,\ldots,D$ are divisible
by $p$. For this case therefore, indeed for the cases of this $y$, in which
the terms of this series are not divisible by $p$,
by necessity the other factor $ax^n-y^n$ will not be divisible by $p$.

\begin{center}
{\Large Corollary 1.}
\end{center}

73. Therefore for whatever number is taken for $x$, if only
not divisible by $p$, for $y$ is always given a value $<p$ which renders
the formula $ax^n-y^n$ divisible by $p$. And in a similar way, if any number
we please is taken for $y$, it can be demonstrated that for $x$ a number $<p$
can always be found by which that formula avoids divisibility by $p$.

\begin{center}
{\Large Corollary 2.}
\end{center}

74. Therefore if $a^m-1$ would be divisible by a prime number $mn+1=p$,
and for $x$ were taken any number $b$ not divisible by $p$, a number
$y$ can always be found so that the form $ab^n-y^n$, or $y^n-ab^n$,
shall be divisible by $p=mn+1$.

\begin{center}
{\Large Corollary 3.}
\end{center}

75. In a similar way if the form $a^m-1$ were to be divisible by a prime number
$p=mn+1$, and for $y$ were taken any number $c$ not divisible by $p$, a
number $x$ will always be able to be found so that the form $ax^n-c^n$,
or $c^n-ax^n$, shall be divisible by $p=mn+1$.

\begin{center}
{\Large Theorem 20.}
\end{center}

76. If the form $ab^n-c^n$, or $c^n-ab^n$, were to be divisible by a prime
number $p=mn+1$, then by taking $d$ to be any number we please, providing
that it were not divisible by $p$, a number $x$ will always be able to be found
so that either the form $ax^n-d^n$, or $ad^n-x^n$, or $d^n-ax^n$,
or $x^n-ad^n$ shall be divisible by this prime number $p=mn+1$.

\begin{center}
{\Large Demonstration.}
\end{center}

Since this form $ab^n-c^n$, or $c^n-ab^n$, was divisible by a prime number
$p=mn+1$, then further the number $a^m-1$ will be divisible by the same
prime number $p=mn+1$. (71) Truly if
$a^m-1$ is divisible by $p$, the number $x$ will be given
by taking any number $d$ not divisible by $p$, so that either the form
$ax^n-d^n$, or on the other hand the form $ad^n-x^n$, or $d^n-ax^n$, or
$x^n-ad^n$, shall be thus divisible by the prime number $p=mn+1$.

\begin{center}
{\Large Corollary.}
\end{center}

77. Therefore by putting $d=1$, if the formula $ab^n-c^n$ were divisible by a
prime number $p=mn+1$, then a number $x$ will be given so that either the form
$ax^n-1$, $a-x^n$, or $x^n-a$ shall be divisible by this prime number $p$.

\begin{center}
{\Large Scholion.}
\end{center}

78. The nineteenth theorem, which is the inverse of the eighteenth theorem,
I shall now 
propose in another way, but without demonstration, and even though I have
attempted the
demonstration of it in many ways, I 
had not however been able to uncover it, until I took up the use of this
method:
which therefore is seen as most noteworthy, and there should be no doubt that
this shall be applied to many other arcane journeys of numbers. For this
method, which is encountered in the consideration of differences, was recently
brought to life by me, which I demonstrated in the service of the
most beautiful theorem of Fermat, 
by which every prime number of the form $4n+1$ is affirmed to be the sum of two
squares; before which I had been unable to persevere in any other way.

\end{document}